\documentclass[12pt,leqno]{article}

\usepackage{amsmath, amsthm, amsfonts, amssymb, color}
\usepackage[matrix,arrow]{xy}
\usepackage{mathrsfs}
\usepackage{enumerate}
\newcommand{\E}{\mathbb E}
\newcommand{\F}{\mathcal F}
\newcommand{\R}{\mathbb R}
\newcommand{\N}{\mathbb N}
\newcommand{\Hi}{\mathcal H}

\newcommand{\B}{\mathcal B}

\newcommand*{\pairing}[3]{\sideset{_{#1}^{}}{_{#3}^{}}{\mathop{\left\langle#2\right\rangle}}}

\newcounter{ncounter}[section]

\makeatletter \@addtoreset{equation}{section}\makeatother
\newtheorem{rem}[ncounter]{Remark}
\newtheorem{example}[ncounter]{Example}
\newtheorem{prop}[ncounter]{Proposition}
\newtheorem{defi}[ncounter]{Definition}
\newtheorem{theorem}[ncounter]{Theorem}
\renewcommand{\P}{\mathbb P}

\title{{\bf A Note on variational solutions to SPDE perturbed by Gaussian noise in a general class}}
\author{\bf Michael R\"ockner \footnote{This work has been supported by NSF-Grant 0603742} \\
\footnotesize 
Department of Mathematics and BiBoS, Bielefeld University, Bielefeld, Germany;\\
\footnotesize
Email: roeckner@math.uni-bielefeld.de\\
\and \bf Yi Wang$^{\; \; *} $ \\
\footnotesize
Department of Mathematics, Purdue University, West Lafayette, USA;}

\begin{document}
\maketitle
\begin{abstract}

This note deals with existence and uniqueness of (variational) solutions to the following type of stochastic partial differential equations on a Hilbert space $ \Hi $\\
\begin{equation*}
 dX(t) = A(t,X(t))dt + B(t,X(t))dW(t) + h(t) \; dG(t)
\end{equation*}

\noindent where $ A $ and $ B $ are random nonlinear operators satisfying monotonicity conditions and $ G $ is an infinite dimensional Gaussian process adapted to the same filtration as the cylindrical Wiener pocess $ W(t),t \geq 0 $.

\end{abstract}

\section{Introduction}

Recently there has been an enormous interest in the study of Hilbert space-valued SPDE with fractional Brownian motion (fBm) representing the perturbing noise. The first paper on this topic we are aware of is the work by W.Grecksch and V.V. Anh. \cite{WGrecksch}, followed by many others (see e.g. \cite{du-ma-pa}, \cite{du-ma-pa2}, \cite{du-ma-pa3}, \cite{ma-nu}, \cite{nu-vu}, \cite{nu}, \cite{ti-tu-vi} and the references therein). The main thrust of the research efforts is now to allow the coefficient in the integral with respect to fBm to depend on the solution. But since It\^o stochastic integration does not apply (because fBm is not a semimartingale), one has to use Skorohod stochastic integration or a pathwise integration approach. Both require more restrictive conditions on the coefficients of the equation. The purpose of this paper is another direction, namely sticking with additive noise we generalize the assumptions on the coefficients. In particular, we want to include nonlinear equations which are not semilinear in order to include more general classical PDE perturbed by Gaussian noise. We do this on the basis of \cite{WGrecksch} (i.e. we use the well-known technique of substracting the additive noise), but allow much more general drift coefficients and even start with an SPDE with a standard Wiener type noise which we disturb by another Gaussian (non-semimartingale) type noise. More precisely, our equation looks as follows

\begin{equation}
\label{eq:1.1}
 dX(t) = A(t,X(t))dt + B(t,X(t))dW(t)+h(t)dG(t)
\end{equation}

\noindent with $ A $ and $ B $ being random and satisfying standard monotonicity conditions (cf. (H1)-(H4) below), $ W(t), t \geq 0 $, a cylindrical Wiener process and $ G(t), t \geq 0 $, an infinite dimensional Gaussian process with suitable integrand $ h $ to be specified below.\\ \\
One main contribution in this note is to specify a class of $ G $ and $ h $ as above so that one can use results from the variational approach to solve SDE (\ref{eq:1.1}) by the said technique of subtracting the additional noise

\begin{equation}
\label{eq:1.2}
 \int_0^\cdot h(t) \; dG(t) \text{.}
\end{equation}

Another contribution is to provide new examples of SPDE to which our results apply.\\ \\
Below we only summarize the main results and sketch the ideas of the proofs. Detailed proofs and additional results are contained in \cite{YWang}

\section{A class of Gaussian integrals}

Let ($ \Omega, \F, \P $) be a complete probability space and for $ T > 0 $ fixed let $ g = (g(t))_{t \in [0,T]} $ be a Gaussian process on ($ \Omega, \F, \P $) with covariance function

\begin{equation}
\label{eq:2.1}
 R(t,s) := \int^s_0 \int^t_0 \phi(u,v) \; du \; dv ; \quad t,s \in [0,T] \text{,}
\end{equation}

where $ \phi $ satisfies the following condition:\\

$ (C_R) \phi \in L^1 ([0,T]^2, \; ds \; dt; \; \R), \phi $ is symmetric, positive definite such that there exist $ p \in (1, \infty) $ and $ C \in (0, \infty) $ such that

\begin{equation}
\label{eq:2.2}
 \int^T_0 \int^T_0 |f(u) f(v) \; \phi (u,v)| \;  du \; dv \leq C \| f \|^2_{L^p ([0;T];\R)} \quad \text{ for all } f \in L^p ([0,T]; \R) \text{.}
\end{equation}

\begin{rem}\label{2.1}

One can show that if $ \phi (u,v) = \Psi (u-v), u, v \in [0,T] $ for some $ \Psi \in L^r ([0,T]; \R) $ with $ r \in (1, \infty) $. Then the inequality in condition $ (C_R) $ holds with $ p:= \frac{2r}{2r-1} $.
In particular, $ (C_R) $ holds for fBM with Hurst parameter $ H \in (\frac{1}{2}, 1) $ where $ \Psi (u) = H(2H-1) |u|^{2H-2} $ (see \cite[Lemma 2.0.2 and Example 2.0.3]{YWang}).

\end{rem}

Under Condition $ (C_R) $ for any separable real Hilbert space $ (\Hi, \langle \; , \; \rangle_{\Hi}) $ and $ f \in L^p ([0,T]; \Hi) $ in the usual way one can define a stochastic integral $ Int(f)(t) := \int^t_0 f \; dg \in L^2 (\Omega ; \; \Hi), \; t \in [0,T]$, such that for all $ f,h \in L^p ([0,T]; \; \Hi) $

\begin{equation}\label{eq:2.3}
 \E \left \langle \int^T_0 f(s) \; dg(s), \; \int^T_0 h(s) \; dg(s) \right \rangle_{\Hi} = \int^T_0 \int^T_0 \langle f(s), h(s') \rangle_{\Hi} \; \phi (s,s') \; ds \; ds'\text{.}
\end{equation}

$ Int(f) $ has $ \P $-a.s. continuous sample paths in $ \Hi $, provided $ f \in L^{p+\epsilon} ([0,T]; \; \Hi) $  for some $ \epsilon > 0 $ (see \cite[Propositions 2.4 and 2.8]{YWang}).
\\
\\
Now let us pass to an infinite dimensional Gaussian process $ G (\text{replacing}\; g) $ taking values in a separable Hilbert space ($ U, \langle \; , \; \rangle_U $) by taking independent copies $ g_n, n \in \N $, of g on a complete probability space ($ \Omega , \F , \P $), fixing an orthonormal basis $ e_n, n \in \N  $, of U and defining

\begin{equation}\label{eq:2.4}
 G(t) := \sum^{\infty}_{n=1} \sqrt{\lambda_n} g_n (t) e_n , t \in [0,T],
\end{equation}

where $ \lambda_n > 0 $ are such that

\begin{equation}\label{eq:2.5}
 \sum^{\infty}_{n=1} \sqrt{\lambda_n} < \infty.
\end{equation}

We note that as usual when dealing with $ g_n $ with covariance given by (\ref{eq:2.1}) (i.e. $ g_n $ is not a Browniem motion) one needs the stronger condition (\ref{eq:2.5}) rather than just assuming \[ {\sum\limits^{\infty}_{n=1}} \lambda_n < \infty \text{.} \]\\
Then for $ h \in L^p([0,T]; L(U,\Hi)) $, where $ L(U, \Hi) $ denotes the set of all bounded linear operators from $ U $ to $  \Hi $, we can define

\begin{equation}\label{eq:2.6}
 \int^t_0 h(s) \; dG(s) := \sum^{\infty}_{n=1} \sqrt{\lambda_n} \int^t_0 h(s) e_j \; dg_j(s), \quad t \in [0,T],
\end{equation}

and one proves (cf. \cite[Lemmas 2.0.10, 2.0.11 and Corollary 2.0.13]{YWang})

\begin{prop}\label{2.2}

$ \int^\cdot_0 h(s) \; dG (s) $ is an $ \Hi $-valued Gaussian process, which, if\\
$ h \in L^{p+\epsilon}([0,T]; \; L(U,\Hi)) $ for some $ \epsilon > 0 $, has $ \P$-a.s. continuous sample paths. Furthermore, for all $ x,y \in \Hi ; h_1 , h_2 \in L^{p}([0,T]; \; L(U, \Hi)) $

\begin{align*}
 & \E \left ( \left \langle \int^T_0 h_1(s) dG(s), x \right \rangle_{\Hi} \quad \left \langle \int^T_0 h_2 (s) \; dG(s), y  \right \rangle_{\Hi} \right )
\\
 = & \int^T_0 \int^T_0 \left \langle h_2(s') Q h_1(s)^* x,y \right \rangle_{\Hi} \phi (s,s') \; ds \; ds'
\end{align*}

where  $ Q z := {\sum\limits_{n=1}^{\infty}} \lambda_n \langle z, e_n \rangle_U \; e_n , \; z \in U, $ is the covariance operator of $ G $. Finally,

\begin{align*}
 & \E \left \langle \int^T_0 h_1(s) dG(s), \quad \int^T_0 h_2 (s) dG(s) \right \rangle_{\Hi}
\\
 = & \int^T_0 \int^T_0 Tr ( h_2(s') Q h_1(s)^*) \phi (s,s') \; ds \; ds' .
\end{align*}

\end{prop}

\section{Main result}

Let ($ H, <, >_H $) be a separable real Hilbert space with dual $ H^* $ and $ V $ a reflexive Banach space such that $ V \subset H $ continuously and densely. Then for its dual $ V^* $ we have

\begin{equation}\label{3.1}
 V \subset H \subset V^*
\end{equation}

continuously and densely, i.e. ($ V, H, V^* $) forms a Gelfand triple. We note that then $ V^* $ and hence $ V $ are also separable. Let ($ U, <, >_U $) be another separable real Hilbert space and $ (W(t))_{t \in [0,T]} $ a cylindrical Wiener process on $ U $ on some stochastic basis ($ \Omega , \F , (\F_t)_{t \in [0,T]}, \P $). Let $ G $  be defined as in Section 2 on ($ \Omega , \F , \P $) (in particular, $ C_R $ is assumed to hold) and suppose that $ G $ is ($\F_t$)-adapted.\\
\\
Now let us fix the conditions on the coefficients $ A , B , h $ in (\ref{eq:1.1}). Let $ L_2(U,H) $ denote the set of all Hilbert-Schmidt operators from $ U $ to $ H $ and let 

\begin{equation*}
 A:[0,T] \times V \times \Omega \rightarrow V^* ,\; B:[0,T] \times V \times \Omega \rightarrow L_2 (U,H)
\end{equation*}

\noindent be \emph{progressively measurable}, i.e. for every $ t \in [0,T] $, these maps 
restricted to $ [0,t] \times V \times \Omega $ are $ \B ([0,t]) \otimes \B (V) \otimes \F_t$-measurable. 
As usual by writing $A(t,v)$ we mean the map $ \omega \mapsto A(t,v,\omega)$.
Analogously for $B(t,v)$. We impose the following conditions on $A$ and $B$:

\begin{enumerate}

\item [(H1)] (\emph{Hemicontinuity}) \index{hemicontinuity} For all $u,v,x \in V,\, \omega \in \Omega $ and $t \in [0,T]$ 
   the map 
   \[
      \R \ni \lambda \mapsto \pairing{V^*}{A(t,u+ \lambda v,\omega),x}{V}
   \]
   is continuous.

\item [(H2)] (\emph{Weak monotonicity}) \index{weak monotonicity} There exists $c\in \R$ such that for all $u,v\in V$

\begin{align*}
      & 2 \pairing{V^*}{A(\cdot, u) - A(\cdot ,v) , u-v}{V}
      + \,\| B(\cdot, u) - B(\cdot,v)\|^2_{L_2(U,H)}\\
      & \leq c\|u-v\|_H^2 \text{ on }[0,T]\times \Omega
\end{align*}

\item [(H3)] (\emph{Coercivity}) \index{coercivity} There exist $\alpha \in \,(1,\infty), \, c_1 \in \R,\,c_2\in\, (0,\infty)$ 
   and an $(\F_t)$-adapted process $f \in L^1([0,T]\times \Omega, dt \otimes P)$ such 
   that for all $v \in V, t \in [0,T]$
   \[
      2 \pairing {V^*}{A(t,v),v}{V}+ \| B(t,v)\|_{L_2(U,H)}^2 \leq c_1 \| v \|_H^2 - c_2 \| v \|_V^{\alpha} + f(t)\quad{on }\;\Omega.
   \]

\item [(H4)] (\emph{Boundedness}) \index{boundedness} There exist $c_3 \in [0,\infty)$ and an $(\F_t)$-adapted 
process
   $ g \in L^{\frac{\alpha}{\alpha -1}}([0,T]\times \Omega , dt \otimes P)$ such that for all 
   $v \in V,\, t\in [0,T]$
   \[
      \|A(t,v)\| _{V^*}\leq g(t) + c_3 \| v\| _V^{\alpha-1} \quad\text{on }\, \Omega,
   \]
   where $ \alpha $ is as in (H3).

\item [(H5)] $ h \in L^{p+\epsilon}([0,T]; L(U, \Hi )) $ for some $ \epsilon > 0 $, where $ \Hi := V $ if $ V $ is a Hilbert space or $ \Hi := H $ if $ h = P.h $ for some orthogonal projection in $ H $ with finite dimensional range in $ V $.

\end{enumerate}

\begin{defi}
\label{def:3.1}
A continuous $ H $-valued ($ \F_t $)-adapted process $ X=(X(t))_{t \in [0,T]}$
is called a (variational) \emph{solution to} (\ref{eq:1.1}), if 

\begin{enumerate}
 \item[(i)] $ \hat X \in L^\alpha ([0,T]\times \Omega, dt \times \P;V) \cap L^2 ([0,T]\times \Omega, dt \times \P; H) $
\end{enumerate}

where $ \alpha $ is as in (H3) and $ \hat X $ is a $ dt \times \P $-equivalence class of $ X $ and $ \P $-a.s.

\begin{enumerate}
 \item[(ii)] $ X(t) = X(0) + \int_0^t A(s,\bar X(s)) ds + \int_0^t B(s, \bar X(s)) dW(s) + \int^t_0 h(s) dG(s) , $
\end{enumerate}

where $\bar X$ is any $V$-valued progressively measurable $ dt \otimes \P $-version of $\hat X$.
\end{defi}

\noindent Now we claim and prove the main result:

\begin{theorem}
 \label{theorem:3.2}
Under the above setting, for any given $ X_0 \in L^2 (\Omega , \F_0 , P, H )$, there exists a unique (variational) solution to (\ref{eq:1.1})
\end{theorem}

\begin{proof}[Proof (sketch).] Let $ w(t), t \in [0,T] $, denote (the continuous version of $ \int^\cdot_0 h(s) \; dG(s) $ and define 

\begin{equation*}
 \bar{A} : [0,T] \times V \times \Omega \rightarrow V^* , \; \bar{B} : [0,T] \times V \times \Omega \rightarrow L_2 (U,H)
\end{equation*}

\noindent by

\begin{align*}
 \bar{A} (t, v, \omega) := & A(t,v+w (t)(\omega),\omega)\\
 \bar{B} (t, v, \omega) := & B(t,v+w (t)(\omega),\omega) .
\end{align*}

Then $ \bar{A}, \bar{B} $ are again progressively measurable (see \cite[Lemma 3.4]{YWang}) and, since $ w(t), \; t \in [0,T] $, has strong moments of all orders as a Gaussian process, one can check that $ \bar{A}, \bar{B} $ also satisfy (H1)-(H4) (cf. \cite[Section 3]{YWang}). Furthermore, $ X $ is a solution to (\ref{eq:1.1}) if and only if

\begin{equation*}
 Y(t) := X(t) - w(t), \; t \in [0,T]
\end{equation*}

\noindent is a solution to

\begin{equation}
\label{3.2}
 d Y(t) = \bar{A} (t, Y(t)) dt + \bar{B} (t, Y(t)) dW(t) .
\end{equation}

But since $ \bar{A}, \bar{B} $ satisfy (H1)-(H4), due to a general theorem on (variational) solutions to SPDE (see e.g. \cite[Theorem 4.2.4]{pr-roe} and the original paper \cite{kr-ro}) there exists a unique (variational) solution to (\ref{3.2}) (see \cite[Section 3]{YWang} for details).

\end{proof}

\begin{rem}
\label{3.3}

It is possible to derive explicit bounds on the second strong moment of the solution in Theorem 2 (see \cite[Section 4]{YWang}).

\end{rem}

\begin{example}
\label{3.4}
Since A and B in (\ref{eq:1.1}) satisfy the standard conditions (H1)-(H4) there are plenty of examples known from classical results in PDE (cf. e.g. \cite[Subsection 4.1]{pr-roe}) and more recent results on SPDE (see e.g. \cite{bo-roe}). We refer to Section 5 in \cite{YWang} where a number of them are worked out. We only mention here that nonlinear operators A of the following type are included:

\begin{align*}
 A(u) = & \triangle \Psi(u)\\
 A(u) = & \text{div} \beta (\triangledown u), 
\end{align*}

\noindent with mappings $ \Psi : \R \rightarrow \R , \; \beta : \R^d \rightarrow \R^d $ respectively, which are continuous, polynomially bounded, monotone and coercive.

\end{example}

\end{document}